\documentclass{amsart}
\usepackage{url}
\usepackage{breakurl}
\usepackage{comment}
\def\GF{\mathrm{GF}}

\begin{document}
\title{Twelve New Primitive Binary Trinomials}
\author[R.\ P.\ Brent]{Richard P.\ Brent}
\address{Australian National University, Canberra, Australia}
\email{trinomials@rpbrent.com}
\author[P.\ Zimmermann]{Paul Zimmermann}
\address{Inria Nancy - Grand Est, Villers-l\`es-Nancy, France}
\email{Paul.Zimmermann@inria.fr}
\date{}


\thispagestyle{empty}

\begin{abstract}
We exhibit twelve new primitive
trinomials over $\GF(2)$
of record degrees $42\,643\,801$, $43\,112\,609$,
and $74\,207\,281$.
In addition we report the first Mersenne exponent not ruled
out by Swan's theorem \cite{Swan62} --- namely $57\,885\,161$ ---
for which no primitive trinomial exists.
This completes the search for the currently known Mersenne prime exponents.
\end{abstract}

\maketitle

Primitive trinomials $x^r + x^s + 1$ of degree 
$r \le 32\,582\,657$ were reported in~\cite{BrZi08}.
We have completed a search for all larger Mersenne
prime exponents found by the GIMPS project~\cite{GIMPS}.
Twelve new primitive trinomials were found
(see Table~\ref{tab:primitive}).
\begin{table}[h]
\begin{center}
\begin{small}
\begin{tabular}{|r|c|l|} \hline
\multicolumn{1}{|c|}{$r$}
                & $s$                   & \multicolumn{1}{c|}{Date} \\ \hline
        42\,643\,801 & 55981, 3706066, 3896488, 12899278, 20150445
                                        &2009 \\ 
        43\,112\,609 & 3569337, 4463337, 17212521, 21078848
                                        &2009 \\ 
        57\,885\,161 & none
                                        &2013 \\ 
        74\,207\,281 & 9156813, 9999621, 30684570
                                        &2016 \\ 
\hline
\end{tabular}
~\\[5pt]
\caption{New primitive trinomials $x^r + x^s + 1$ of degree
a Mersenne exponent $r$, for $s \le r/2$.
For smaller exponents,
see references in~\cite{BrZi08} or our web site~\cite{BrZi-web}.
}
\label{tab:primitive}
\end{small}
\end{center}
\end{table}

\vspace*{-5mm}
Our search used the algorithm of~\cite{BrZi07a},
relying on fast arithmetic in ${\rm GF}(2)[x]$; details
are given in~\cite{BrGaThZi08}.
For the squaring of polynomials over ${\rm GF}(2)[x]$,
we used (since 2016) the new
\verb|_pdep_u64| Intel intrinsic, which gave a speedup of a factor
about 2.5 over the algorithm described in~\cite[\textsection 4]{BrLaZi03}.
On a 3.3Ghz Intel Core i5-4590,
together with improvements in the \texttt{gf2x} library, we were able to
square a degree-$74\,207\,280$ polynomial in about 2 milliseconds,
and to multiply two such polynomials in about 700 milliseconds.
As in \cite{BrZi08}, we produced certificates for non-primitive trinomials
(a certificate is simply an encoding of a nontrivial factor of smallest
degree). The certificates were checked independently with Magma and NTL.
A 3.3Ghz Intel Core i5-4590 takes only 22 minutes to check the certificates of
all $37\,103\,637$ reducible trinomials ($s \le r/2$) of degree $r=74\,207\,281$
with our \texttt{check-ntl} program based on NTL~\cite{NTL980},
the largest factor having degree $19\,865\,299$ for $s = 9\,788\,851$.

\par
\noindent
\textsc{Acknowledgements.}
The authors thank Allan Steel, who independently verified with Magma
the twelve new primitive trinomials given in Table~\ref{tab:primitive},
and Gr\'egoire Lecerf, who verified with Mathemagix the three primitive
trinomials of degree $74\,207\,281$.
Part of the computations reported in this paper were carried out using a
cluster funded by the French ANR CatRel.
Dan Bernstein, Bill Hart, Alex Kruppa and Tanja Lange contributed some
CPU cycles towards the searches for degrees $43\,112\,609$  
and $57\,885\,161$. 
The work of the first author was supported in part by
Australian Research Council grant DP140101417.

\providecommand{\bysame}{\leavevmode\hbox to3em{\hrulefill}\thinspace}
\providecommand{\MR}{\relax\ifhmode\unskip\space\fi MR }
\providecommand{\MRhref}[2]{%
  \href{http://www.ams.org/mathscinet-getitem?mr=#1}{#2}
}
\providecommand{\href}[2]{#2}


\begin{thebibliography}{10}

\bibitem{BrZi-web}
Richard~P. Brent, \emph{Search for primitive trinomials (mod 2)},
  \url{http://maths-people.anu.edu.au/~brent/trinom.html}, 2008.

\bibitem{BrGaThZi08} 
Richard Brent, Pierrick Gaudry, Emmanuel Thom\'e, and Paul Zimmermann,
  \emph{Faster multiplication in ${\rm GF}(2)[x]$}, Proc.~of the 8th
  International Symposium on Algorithmic Number Theory (ANTS {VIII}),
  \emph{Lecture Notes in Computer Science} \textbf{5011}, 
  Springer-Verlag, 2008, 153--166.

\bibitem{BrLaZi03} 
Richard~P. Brent, Samuli Larvala, and Paul Zimmermann, \emph{A fast algorithm
  for testing reducibility of trinomials mod $2$ and some new primitive
  trinomials of degree $3021377$}, Math.\ Comp.\ \textbf{72} (2003), 
  1443--1452.
  MR1972745 (2004b:11161) 

\bibitem{BrZi07a} 
Richard Brent and Paul Zimmermann,
  \emph{A multi-level blocking distinct
  degree factorization algorithm}, Contemporary Mathematics
  \textbf{461} (2008), 47--58.

\bibitem{BrZi08}
\bysame,
  \emph{Ten new primitive trinomials},
  Math.\ Comp.\ \textbf{78} (2008), 
  1197--1199.

\bibitem{BrZi11}
\bysame,
  \emph{The great trinomial hunt},
  Notices Amer.\ Math.\ Soc.\ \textbf{58:2} (2011), 233--239.

\bibitem{GIMPS}
The Great Internet Mersenne Prime Search, 
\url{http://www.mersenne.org}.

\bibitem{KuMa91}
Y.~Kurita and M.~Matsumoto, \emph{Primitive $t$-nomials $(t = 3, 5)$ over ${\rm
  GF}(2)$ whose degree is a {M}ersenne exponent $\le 44497$}, 
  Math.\ Comp.\ \textbf{56} (1991), 817--821.
  MR1068813 (91h:11138)

\bibitem{NTL980}
Victor Shoup, \emph{{NTL}: A library for doing number theory},
  \url{http://www.shoup.net/ntl/}, 2016.

\bibitem{Swan62}
R.~G. Swan, \emph{Factorization of polynomials over finite fields}, Pacific J.
  Math. \textbf{12} (1962), 1099--1106.
  MR0144891 (26 \#2432)

\end{thebibliography}
\end{document}